**Decoding gender differences: Intellectual profiles of children with specific learning disabilities**

Authors: D. Giofrè[1*†], K. Allen[2], E. Toffalini[3], I. C. Mammarella[4], & S. Caviola[4]

David Giofrè  http://orcid.org/0000-0003-1145-8642

Katie Allen  http://orcid.org/0000-0003-1145-8642

Enrico Toffalini  http://orcid.org/0000-0002-1404-5133

Irene Cristina Mammarella  http://orcid.org/0000-0002-6986-4793

Sara Caviola  http://orcid.org/0000-0002-4556-3179

**Affiliations:**

[1]Disfor, University of Genoa, Italy

[2]Department of Psychology, University of Durham

[3]Department of General Psychology, University of Padua, Italy

[4]Department of Developmental and Social Psychology, University of Padova

*David Giofrè, Disfor, University of Genoa, Italy, Corso Andrea Podestà, 2, david.giofre@gmail.com

**Decoding gender differences: Intellectual profiles of children with specific learning disabilities**

**Abstract**

There has been a significant amount of debate around gender differences in intellectual functioning, however, most of this research concerns typically developing populations and lacks research into atypically developing populations and those with specific learning disabilities (SLD). To address this, we examined performance on the WISC-IV in children with SLDs (N=1238, N female= 539, Age range = 7-16 years). We further divided the sample into those with specific deficits in reading, mathematics, and those with mixed disorder. Results indicate that gender predicts significant differences in the working memory index and processing speed index only, indicating a small but significant female superiority. Results also show different profiles for the different disorders investigated, with some gender differences emerging. The most prominent gender difference appears to be in the coding subtest indicating a female advantage, particularly in those with SLDs with mathematical difficulties. We discuss the theoretical and clinical implications of the findings.

**Decoding gender differences: Intellectual profiles of children with specific learning disabilities**

There has been a long-lasting debate on the possible presence of gender differences in general intellectual functioning. Research in this area has traditionally used intelligence tests, particularly Wechsler scales for children (WISC) and adults (WAIS), to test for the possible presence of differences in general intellectual functioning, or *g*, and other subtests. Most of these studies have been performed on standardization samples from several different countries. However, few studies have also taken into consideration other populations, including children with various disabilities and in particular children with specific learning disabilities (SLDs).

Several different accounts have been proposed to explain gender differences in intelligence and more broadly in cognitive abilities. Influential theories proposed that cognitive differences between genders, particularly in spatial tasks, are probably due to differences in brain lateralization and hormone levels, however this evidence seems to be rather inconclusive (see Miller & Halpern, 2014 for a review). Some evidence suggests that androgens might affect spatial ability in females, not only directly but also indirectly, through male-typed activity interests (Berenbaum et al., 2012). Absolute and relative brain size seems to be an important determinant of intelligence (Roth & Dicke, 2005), and it has also been advocated as a possible explanation of differences in general intelligence (Jensen & Johnson, 1994). Biological theories, however, have been sharply criticized by other influential authors claiming that group differences in IQ probably have an environmental origin (Nisbett et al., 2012). While a comprehensive review of explanatory theories on gender differences in general intelligence is outside the scope of this paper (see Halpern & Wai, 2019 for a recent review), it is worth mentioning that none of these theories are conclusive and the question of whether men and women are fundamentally different or similar requires further investigation (see Hyde 2014).

Concerning differences in the general intelligence quotient, evidence using the WAIS-IV, WAIS-III and WAIS-R generally shows a small difference in full scale IQ (FSIQ) favoring males, in particular in China, USA, Canada, Japan and Italy (Dai et al., 1991; Hattori & Lynn, 1997;

Longman et al., 2007; Pezzuti et al., 2020; see Lynn, 2017 for a review). An empirical investigation of gender differences using several different tasks and batteries (including the WAIS) concluded that it is probable that differences are small on the general factor but larger on other secondary factors (i.e., verbal, perceptual and mental rotation factors), with males performing better on mental rotation tasks while females display better memory (Johnson & Bouchard, 2007). Looking at specific subtests – over and above the effect of the general factor – large differences were found in several tests, including, for example, block design and coding (Johnson & Bouchard, 2007).

Results looking at children's performance on the WISC show a different pattern, with some countries showing more pronounced differences, while others show only trivial or negligible differences (see Chen & Lynn, 2020 for a discussion). Past evidence using the WISC-R showed that, counterintuitively, males tended to outperform females on verbal abilities, particularly on the information subtest (Born & Lynn, 1994; Lynn et al., 2005; see Lynn, 2021 for a review). Recent evidence using the WISC-IV with Italian children shows a similar trend with a male advantage on most verbal tasks (Pezzuti & Orsini, 2016). Other results indicate that males perform better than females on the block design and object assembly subtests of the WISC-R (e.g., Lynn et al., 2005). Evidence for memory tasks is more mixed, with some authors finding a female advantage on some verbal memory tasks (Jensen & Reynolds, 1983), while others find no differences whatsoever (Pezzuti & Orsini, 2016; see Roivainen, 2011 for a review).

Finally, concerning processing speed tasks, girls tend to consistently outperform boys in coding, and this is true across different samples and different versions of the WISC (Born & Lynn, 1994; Goldbeck et al., 2010; Jensen & Reynolds, 1983; Lynn & Mulhern, 1991; Pezzuti & Orsini, 2016). The reason behind this advantage, which has been consistently observed, is rather elusive, but some authors claim that this generally reflects faster processing in writing speed and associated learning, and faster retrieval from secondary memory (e.g., Lyle & Johnson, 1974; see also Halpern & Wai, 2019). To sum up, recent data on FSIQ using the WISC-IV in children with typical development showed either no difference or a very small difference in IQ favoring males, while some differences remained at the level of the indices and subtests but not on general intelligence (e.g., Goldbeck et al., 2010; Pezzuti & Orsini, 2016).

SLDs are neurodevelopmental disorders that lead to persistent difficulties in several academic domains, in particular reading, spelling, and calculation. Beyond impairments in the academic domains, children with SLD also present with a series of other deficits in other cognitive domains, including working memory and processing speed, which are the object of the assessment of several intelligence batteries (see Cornoldi & Giofrè, 2014 for a review). The consideration of intelligence has always been crucial for the diagnosis of specific learning disorders, and is routinely included in the diagnostic process of children with SLD.

Several different studies confirm the presence of an atypical cognitive profile in children with SLD. According to the DSM-5, SLDs are more frequently diagnosed in males than in females (American Psychiatric Association, 2013). However, data on the general population with extremely large samples indicate that there is a trade-off with more males having difficulties with reading and more females presenting with difficulties in mathematics (Giofrè et al., 2020; Stoet & Geary, 2013). Looking at different diagnostic categories within SLDs, there seem to be more males than females with SLDs with impairments in reading (Reilly, 2012). Concerning children with SLDs with impairments in mathematics, the prevalence is somewhat less clear, but in general females seem to be overrepresented in this category (e.g., Shalev et al., 2000).

Recent evaluations of the cognitive profile of children with SLD revealed that these children tend to have an atypical profile of strengths and weaknesses. The WISC-IV has been one of the most widely used batteries for the assessment of children with typical development and with SLD (Evers et al., 2012). This particular battery encompasses 10 core subtests with four main indices measuring ability in several distinct domains including verbal (verbal comprehension index, VCI), visuo-perceptual (perceptual reasoning index, PRI), working memory (working memory index, WMI), and processing speed (processing speed index, PSI). Alongside these principal indices, intelligence is traditionally measured through a single measure: the full-scale intellectual score (FSIQ, an indicator of the *g*-factor).

Children with SLD tend to have higher scores on the VCI and PRI and lower scores on WMI and PSI (e.g., Giofrè et al., 2017). This pattern has been repeatedly confirmed in various populations (Cornoldi et al., 2014; De Clercq-Quaegebeur et al., 2010; Poletti et al., 2016). This pattern of results has led to the suggestion that the general ability index (GAI), composed of the VCI and PRI, should be preferred over the FSIQ which also includes WMI and PSI (which together compose an additional ancillary index called the cognitive proficiency index; CPI) (Giofrè et al., 2017; Poletti, 2016).

The proficiency deficits of children with SLD very often relate to working memory (WM) (Swanson & Ashbaker, 2000; Swanson, 1993), and processing speed (PS) (Proctor, 2012). A recent meta-analysis on children with SLD confirms the presence of deficits on various WM tasks, but particularly when tasks involve the numerical and the verbal domains (Peng & Fuchs, 2016). Both working memory and processing aspects seem to be crucially impaired in children with SLD (Swanson, & Sachse-Lee, 2001). A tentative explanation for what may drive these differences has been proposed that, in typically developing children, females tend to outperform males on tasks requiring rapid access to, and use of, information from long-term memory (see Halpern & Wai, 2019 for a discussion). Results using the WISC-R on children with SLD align closely with this explanation, showing that girls have a superior performance on WM (e.g., digit span) and PS (e.g., coding) tasks (e.g., Spafford, 1989; Vance et al., 1980).

It is also noteworthy that SLD is a very broad category, encompassing children with different difficulties, including difficulties with reading, with mathematics, or with more than one diagnosis within this category. Despite some similarities between different diagnoses within each SLD category, some differences also emerge. Recent results focusing on performance on the WISC-IV of children with various diagnoses within the SLD category showed that performance on this battery has some similarities (e.g., the GAI generally higher than the FSIQ, and the CPI), but also differs to some extent between different diagnoses (Toffalini, Giofrè, & Cornoldi, 2017). In their report, Toffalini and co-authors showed that children with specific reading disabilities seem to be more impaired on verbal aspects (on the majority of tasks included in the VCI). Conversely, children with SLD with mathematics impairments perform somewhat better on verbal tasks but tend to perform worse on visuospatial processing tasks (i.e., the majority of tasks included in the PRI). In addition, children with SLD with mixed disorder, which very often includes children with problems in both mathematics and reading, represent a somewhat middle ground, characterized as having lower intellectual functioning, including a lower FSIQ, as well as a generally lower performance in all indices and subtests. Despite the theoretical and practical relevance of this specific population the performance of males and females on the WISC in children with SLD has seldom been examined.

An important limitation of most of the previous research is that results are based on performance at the observed level. Summary scores, such as FSIQ scores, along with the other principal indices of Wechsler scales, are determined by a weighted combination of scores on subtests and, if males have higher scores on some subtests and females on other subtests, then, depending upon the weights assigned to each subtest, you could produce a summary score that favored men over women (Hunt, 2011). For these reasons, the comparison between measurements for males and females (factor scores) on intelligence (i.e., on $g$), should be performed at the latent level rather than on indicators such as indices derived from Wechsler batteries (see for Pezzuti & Orsini, 2016 and Pezzuti et al, 2020 for a similar procedure). The argument is that the weighting of individual subtests will then be done by rational analysis of the data, rather than by using weights that were arbitrarily assigned to the subtests (see Hunt, 2011). It is also important to note that there are several advantages to measuring latent rather than manifest variables including, but not limited to, i) latent models allow us to obtain more precise estimates and to compensate for some biases due to errors related to variance on manifest variables, ii) differences in factors are weighted depending on the loading of a factor on its manifest variables, thus producing more realistic estimates of the latent ability; and iii) manifest variables are in general less reliable as compared to latent factors (see Kline, 2011 for a more in depth description of the advantages of using latent rather than manifest variables). Another limit of the aforementioned studies is that they are limited to children with typical development and rarely on children with disabilities (e.g., specific learning disability, SLD).

Using the WISC-III on a sample of 440 children with SLD, Slate (1997) found that males exhibited statistically higher full scale, verbal, and performance IQs than did females as well as higher scores on six of seven subtests, with females outperforming males only on the coding subtest. On four subtests, arithmetic, similarities, picture arrangement, and digit span, the mean gender difference was small and not statistically significant. Three subtests were the most important in discriminating between groups: object assembly, coding, and information. This study was one of the few investigating gender differences in a sample of SLD using the WISC, however, a few issues remained unclear. First, the sample included children with SLD, but without distinguishing between different SLD subtypes, and some subtests are no longer in use in later versions of the WISC (e.g., information).

Based on the literature presented so far, we decided to perform a study investigating gender differences on the WISC-IV in children with SLD. With this study we aimed to evaluate standardised differences in performance, both on indices (at the observed and at the latent level), and on subtests. We also aimed to compare the cognitive profile of male and female children with different diagnoses of SLD to verify whether the profiles are similar or differ across genders and across diagnoses.

## Method

**Participants**

Data on a large number of children with SLD was collected under the auspices of the Italian Association for Learning Disabilities (AIRIPA). A subset of this data has been included in previously published articles (Cornoldi et al., 2014; Giofrè et al., 2016, 2017; Giofrè, Pastore, et al., 2019; Giofrè, Toffalini, et al., 2019; Giofrè & Cornoldi, 2015; Toffalini, Giofrè et al., 2017; Toffalini, Pezzuti et al., 2017), however, those manuscripts did not address the issue examined in the present study. All children received a diagnosis within the F81 category (i.e., specific developmental disorders of scholastic skills) of the ICD-10 International Coding System (World Health Organization, 1992), which is the classification system generally consulted in Italy for SLD. This was done following the guidelines indicated by the National Italian Consensus Conference on SLD published by the Italian Ministry of Health (Istituto Superiore di Sanità, 2011). Experts were invited to provide data on children with a diagnosis of SLD, but to exclude cases in which a comorbid neuropsychological diagnosis was also present (e.g., attention-deficit hyperactivity disorder, developmental coordination disorder, or specific language disorder).

We included participants who received a diagnosis of dyslexia, dyscalculia, or mixed disorder according to the ICD-10 coding system (N = 1,236): 455 (177 females) children with a reading disorder (F81.0); 120 (80 females) children with a specific disorder of arithmetical skills (F81.2); and 661 (280 females) children with a mixed disorder of scholastic skills (F81.3). Children with SLD were between 7.1 and 16.9 years of age. Data were provided for 537 females ($M_{age}$ = 11.75, $SD$ = 2.42) and 699 males ($M_{age}$ = 11.66, $SD$ = 2.45). Females and males did not differ statistically in age, Cohen's $d$ = 0.04 (95%CI -0.08, 0.15), $t(1234)$ = 0.64, $p$ = .52.

**Instrument**

The Italian version of the WISC-IV (Orsini et al., 2012) was used. For the purposes of the present study, we examined the scores on the 10 core subtests of the WISC-IV, i.e., block design (BD), similarities (SI), digit span (DS), picture concepts (PC), coding (CD), vocabulary (VC), letter–number sequencing (LN), matrix reasoning (MR), comprehension (CO), and symbol search (SS). We also calculated the Full-Scale IQ (FSIQ), and the four main indices: verbal comprehension index (VCI), perceptual reasoning index (PRI), working memory index (WMI), and processing speed index (PSI). We then calculated the scores for two additional indices, i.e., general ability index (GAI), obtained from VCI and PRI; and the cognitive proficiency index (CPI), obtained from WMI and PSI.

**Data analytic approach**

The R program (R Core Team, 2021) with the "lavaan" library (version 0.5-17; Rosseel, 2012) was used for structural equation modelling (SEM). Model fit was assessed using various indices according to the criteria suggested by Hu and Bentler (1999). We considered the chi-square ($\chi^2$), the comparative fit index (*CFI*), the non-normed fit index (*NNFI*), the standardized root mean square residual (*SRMR*), and the root mean square error of approximation (*REMSEA*) (Kline, 2011).

Analyses of the visual matching tasks were performed using generalized linear mixed models (Pinheiro & Bates, 2000) using the "lm4" package (Bates et al., 2015). GLMM is a robust analysis that allows controlling for the variability of subjects, limiting the loss of information due to by-subject analyses (Baayen et al., 2002). P-values for fixed effects were obtained using the package "car" (Fox & Weisberg, 2019). Figures were obtained using the package "ggplot2" (Wickham, 2016).

To obtain confidence intervals, we employed robust distribution independent bootstrap statistics. According to these, 'significant' differences appear if appropriate 95% bootstrap confidence intervals do not overlap. Hence, the term 'significant' will refer to such differences in confidence intervals. All bootstrap confidence interval estimations used 10,000 permutations with replacement (Chihara & Hesterberg, 2011) and computed bias-corrected and accelerated (BCa) confidence intervals (Efron, 1987). We assessed group differences by computing 95% BCa bootstrap confidence intervals for the main measures, which provide a better statistical solution than simply reporting *p*-values (Cumming, 2013). To interpret standardized effect sizes we used the criteria proposed by Cohen (1988): .20, .50, and .80 for Cohen's *d* were considered to be small, medium, and large effects, respectively. The analyses were performed using the "boot" package (Canty & Ripley, 2019).

## Results

**Descriptive and standardized differences between males and females**

Descriptive statistics and standardized differences (Cohen's *d*) for female vs. male were calculated for all indices and subtests in the three different subgroups (Table 1).

Table 1 about here

**Structural equation modelling (SEM)**

We performed a series of SEM to investigate the effect of gender on latent factors (VCI, PRI, WMI, PSI). In a first model, a traditional four factor model (VCI, PRI, WMI, and PSI) was fitted. This model allowed us to directly test the effect of gender (operationalized as 0 female, 1 male variable) on latent factors (VCI, PRI, WMI, and PSI) (Figure 1). In a second model, we investigated the effect of gender on the *g*-factor. We fitted a model with a first order factor (*g*-factor), and four second order factors (VCI, PRI, WMI, PSI) (Figure 2).

The sample size was constrained by available cases, not determined by a priori calculation. However, we established whether our sample size was sufficient for a minimum effect of interest through power analysis via simulation with 5,000 iterations. Specifically, we simulated the SEM shown in Figures 1 and 2 using the covariance matrix from the Italian normative sample of the WISC-IV as the prior (Orsini et al., 2012). Concerning the minimum effects of interest, we set a coefficient for gender on the latent factors equivalent to a standardized effect of $B = .10$ (corresponding to Cohen's d = .20). As power is constrained by the smallest group, and our female subsample was just over 500, we performed all simulations with N = 1,000. Considering a critical α = .05, power was slightly suboptimal: 82% for the VCI; 73% for the PRI; 69% for the WMI; 75% for the PSI; finally, it was 80% for the g-factor. Different models were simulated for each of the effects. For a slightly larger standardized effect of interest of $B = .15$ (corresponding to Cohen's d = .30), power seemed sufficient: 99% for the VCI; 98% for the PRI; 96% for the WMI; 98% for the PSI; 99% for the g-factor.

*Preliminary analysis of fit of measurement model and invariance*

The fit of the first measurement model (i.e., four-factor) was adequate, $\chi^2(29) = 77.32$, *RMSEA* = .037, *SRMR* = .026, *CFI* = .972, *NNFI* = .968. To assess measurement invariance across the two genders, we fit multigroup models, with subsequent equality constraints of loadings, intercepts, and residuals, across gender. For simplicity, model comparisons were conducted using the Bayesian Information Criterion (BIC, lower is better; Kline, 2011). Configural invariance was demonstrated by the very good fit of the multigroup model (with all parameters freely estimated in the two groups): $\chi^2(58) = 100.58$, *RMSEA* = .034, *SRMR* = .027, *CFI* = .982, *NNFI* = .972, *BIC* = 58,544.02. Metric invariance was established since the model with loadings constrained across groups was largely better, *BIC* = 58,507.13, *ΔBIC* = -36.89. Scalar invariance was not established, however, *BIC* = 58,516.08, *ΔBIC* = +8.95. Nonetheless, partial scalar invariance could be established, as a model with all intercepts constrained across groups, except for Block design (first step) and Coding (second step), reached the best possible fit, *BIC* = 58,481.51, *ΔBIC* = -25.62. A further improvement in fit was reached when residual variances were also constrained across groups, *BIC* = 58,424.38, *ΔBIC* = -57.13. Therefore, measurement invariance was largely established across genders, except for the intercepts of two subtests. Concerning diagnoses, measurement invariance was fully established, with BIC steadily decreasing across all steps: unconstrained model, $\chi^2(87) = 130.79$, *RMSEA* = .035, *SRMR* = .032, *CFI* = .980, *NNFI* = .969, *BIC* = 58,647.54; metric invariance, *BIC* = 58,579.39; scalar invariance, *BIC* = 58,526.27; strict (residuals) invariance, *BIC* = 58,414.03.

Concerning the second measurement model (i.e., higher order), the model fit was also satisfactory, $\chi^2(31) = 93.04$, *RMSEA* = .040, *SRMR* = .033, *CFI* = .974, *NNFI* = .962. Assessment of measurement invariance yielded virtually the same results as for the previous model. Concerning gender, for the unconstrained model, $\chi^2(62) = 116.10$, *RMSEA* = .038, *SRMR* = .032, *CFI* = .977, *NNFI* = .967, *BIC* = 58,531.07; metric invariance, *BIC* = 58,474.60, scalar invariance, except for Block design and Coding, *BIC* = 58,456.02; further constraining of residuals, *BIC* = 58,398.92.

Concerning diagnoses, measurement invariance was fully established: unconstrained model, $\chi^2(93)$ = 154.12, *RMSEA* = .040, *SRMR* = .038, *CFI* = .972, *NNFI* = .959, *BIC* = 58,628.16; metric invariance, *BIC* = 58,523.31; scalar invariance, *BIC* = 58,484.63; strict (residuals) invariance, *BIC* = 58,372.74.

*Effects of gender on latent factors*

Concerning the first model, including gender as a predictor of all latent factors (but not subtests) still showed a satisfactory fit, $\chi^2(35)$ = 127.49, *RMSEA* = .046, *SRMR* = .031, *CFI* = .962, *NNFI* = .940. The model coefficients are shown in Figure 1. Importantly, the effect of gender on VCI and PRI were small and not statistically significant, while gender was predicting a statistically significant, albeit small, portion of the variance on WMI, and on PSI. Based on the portion of explained variance we were able to calculate standardized mean differences on each specific factor: VCI (*d* = 0.113 [-0.019, 0.237], or 1.695 IQ points; *p* = .08), PRI (*d* = -0.056 [-0.203, 0.096], or -0.840 IQ points; *p* = .44), WMI (*d* = -0.183 [-0.319, -0.039], or -2.745 IQ points; *p* = .01), and PSI (*d* = -0.350 [-0.540, -0.136], or -5.250 IQ points; *p* < .001).

Concerning the second model, adding gender as a predictor of the g factor (but not of first-order factors or subtests) still showed acceptable, albeit suboptimal, fit: $\chi^2(40)$ = 173.14, *RMSEA* = .052, *SRMR* = .043, *CFI* = .945, *NNFI* = .925. The model coefficients are shown in Figure 2. Gender was directly predicting the variance of the g-factor, but this effect was not statistically significant and was trivial in terms of magnitude (*d* = -0.074 [-0.223, .0.085], or -1.110 IQ points; *p* = .30).

Figure 1 and 2 about here

**Generalized linear mixed models (GLMM)**

We investigated the effects of gender on performance on the WISC-IV principal indices and subtests considering different SLD profiles. Considering the four main indices (VCI, PRI, WMI and PSI), the effect of index, $\chi^2(3) = 1371.73$, $p < .001$, gender, $\chi^2(1) = 8.85$, $p = .003$, diagnosis, $\chi^2(2) = 118.31$, $p < .001$, index × gender, $\chi^2(3) = 37.42$, $p < .001$, and index × diagnosis, $\chi^2(6) = 36.00$, $p < .001$ were statistically significant, while the effects of gender × diagnosis, $\chi^2(2) = 1.26$, $p = .532$, and index × gender × diagnosis, $\chi^2(6) = 5.87$, $p = .438$, were not statistically significant (Table 1; Figure 3).

We also compared the two ancillary indices (GAI and CPI). The effect of index, $\chi^2(1) = 1178.70$, $p < .001$, gender, $\chi^2(1) = 8.72$, $p = .003$, diagnosis, $\chi^2(2) = 120.58$, $p < .001$, and index × gender, $\chi^2(1) = 23.81$, $p < .001$, were statistically significant, while the effects of index × diagnosis, $\chi^2(2) = 0.55$, $p = .760$, gender × diagnosis, $\chi^2(2) = 1.69$, $p = .429$, and index × gender × diagnosis, $\chi^2(2) = 0.001$, $p = .999$, were not statistically significant (Table 1; Figure 4).

We then compared the performance in the 10 principal subtests (SI, VC, CO, BD, PC, MR, DS, LN, CD, and SS). The effect of index, $\chi^2(9) = 2294.81$, $p < .001$, gender, $\chi^2(1) = 4.98$, $p = .026$, diagnosis, $\chi^2(2) = 113.21$, $p < .001$, index × gender, $\chi^2(9) = 94.83$, $p < .001$, and index × diagnosis, $\chi^2(18) = 79.58$, $p < .001$ were statistically significant, while the effects of gender × diagnosis, $\chi^2(2) = 1.45$, $p = .483$, and index × gender × diagnosis, $\chi^2(18) = 15.09$, $p = .656$, were not statistically significant (Table 1; Figure 5).

Figures 3 to 5 about here

**Gender Differences in WISC-IV Indices Divided by FSIQ Level**

Analyses for different IQ subgroups were also performed and are available online (see Supplemental Online Material).

## Discussion

The main aim of this study was to evaluate the presence of gender differences on the WISC-IV in children with SLD. SEM analyses showed that the difference on the *g*-factor, favoring females, was trivial in terms of magnitude. A similar result was obtained on the indices with the highest *g*-loading: on the VCI, the difference, favoring males, was trivial in terms of magnitude. So was the difference favoring females on the PRI. These results were corroborated both at the observed level and at the latent level confirming that both males and females performed similarly on the g-factor and on the indices with the highest *g*-loading (PRI and VCI). As for WMI and PSI, results showed a female superiority on both indices; a result confirmed both at the observed and at the latent level. Differences were smaller on the WMI and larger on the PSI. Looking at the subtests, males scored higher on four subtests (i.e., SI, VC, CO, which are the indicators of the VCI, and BD), while females outscored males in the remaining six subtests (i.e., PC, MR, DS, LN, CD, and SS). It is worth noting, however, that differences, albeit statistically significant, were generally small in terms of magnitude, apart from CD, in which differences were larger and favoring females.

With some discrepancy across instruments (e.g., Woodcock-Johnson [WJ] cognitive and achievement batteries), the results for processing speed are remarkably stable and robust over time: males scored significantly lower than females across different school-cohorts (Camarata & Woodcock, 2006). One possible explanation can be found in the processing speed definition itself, which can be defined as the ability to automatically perform simple cognitive tasks when under pressure to maintain attention and concentration (Flanagan, McGrew, & Ortiz, 2000). According to this position, males perform worse than females because of their limits in maintaining attention and concentration while performing a simple and repetitive task for an extended period.

The results of gender differences were similar to the results obtained in the Italian standardization sample (Pezzuti & Orsini, 2016). In particular, the directionality of the male/female differences was the same in the two samples, corroborating the evidence obtained with children with typical development. Overall, our results taken together with those in the standardisation sample showed that gender differences, even when statistically significant, tend to not be particularly large. In addition, results obtained in the current report are similar to those obtained in other countries strengthening the idea that gender differences on the *g*-factor are negligible (Colom et al., 2000; Goldbeck et al., 2010; Halpern & LaMay, 2000; Pezzuti & Orsini, 2016; Ranehill et al., 2015; Weiss & Prifitera, 1995). Importantly, and in line with a recent meta-analysis (Roivainen, 2011), differences seemed to be larger on the PSI index, in particular on the coding subtest.

A female superiority on coding has been repeatedly confirmed in almost all standardisation samples with both children and adults, and with children with both typical and atypical development. Notably, the female advantage in coding is confirmed in the WISC-R, WISC-III, WISC-IV, and WISC-V (e.g., Chen & Lynn, 2020; Demo, 1982; Goldbeck et al., 2010; Jensen & Reynolds, 1983; Lynn & Mulhern, 1991; Pezzuti & Orsini, 2016). Lawson and Inglis (1984) suggested that the coding test may be more verbally weighted than previously thought and it is for this reason that the female advantage on this subtest appears (Ryckman, 1981; Vogel, 1990). Other accounts, however, claim that females are generally better on all perceptual speed tasks, not only those with high encodability (Delaney et al., 1981; see also Roivainen, 2011). This latter hypothesis seems to be more consistent with our results and in fact differences on other verbal tasks favour males (e.g., tasks included in the VCI) or favour females but with a small effect size. Some important conclusions can be drawn from these results.

Another possible explanation can be seen by looking at the response dynamics of the specific tasks. In our study, the coding subtest uses visual abstract symbols (also including numbers) that have to be visually inspected and rapidly converted into a series of written elements. Thus, involving fine motor skills, which are important for handwriting. It is worth mentioning that fine and gross motor skills seem to be impaired in children with various SLDs profiles (e.g., Gooch et al., 2015; Rochelle et al., 2006). It is also interesting to note that females seem to outperform males in handwriting skills (e.g., Kono et al., 2019; Tseng & Hsueh, 1997; but see Peters & Campagnaro, 2016 for a different argument). Previous studies accounted for female superiority in writing skills (both speed and legibility) due to their deeper engagement in writing-related activities at school and home (Graham, Berninger, Weintraub, & Schafer, 1998; Lynn & Mikk, 2009). For these reasons, it is possible to speculate that females tend to outperform males in tasks requiring writing-related skills, such as the coding subtest.

Previous evidence indicates that the coding subtest presents with the lowest load on the g-factor in children with SLD, with loadings on g of .14 and .39 in Italian and Irish children with SLD respectively (Giofrè & Cornoldi, 2015; Watkins et al. 2013). Coding, alongside other WISC subtests, was traditionally included because of its power in discriminating between children with typical and atypical development (see Giofrè, Pastore, et al., 2019 for an historical consideration). This observation is also confirmed by that fact that children with various disabilities, including SLD, DCD, autism with lower cognitive abilities and ADHD for example, almost invariably show poor performance on this subtest (Cornoldi et al., 2013; De Clercq-Quaegebeur et al., 2010; Giofrè, Provazza, et al., 2019; Mayes & Calhoun, 2006, 2007; Poletti, 2016; Sumner et al., 2016; Thaler et al., 2013). Looking at the most recent version of the WISC, it is notable that this subtest is still included as one of the core subtests and, having reduced the total number of subtests needed to obtain the IQ score, this test has even more weight, or leverage, on the overall estimate of the FSIQ. The decision to maintain this subtest in the latest version of the WISC might have important consequences for children with various disabilities, particularly for males.

Concerning the second aim of this paper, we aimed to compare the profiles of children with dyslexia, dyscalculia or with a mixed SLD with a particular focus on gender differences. Results of GLMM analyses indicated that main effects, as well as the index by gender interaction, were always statistically significant in all of the models tested. These results confirm that children with different SLD diagnoses tended to have different profiles, that these profiles were not uniform within the various categories, that there were some overall gender differences in the profiles, and that gender differences between subtests and indices were present only in some but not in all subtests. Looking at the graphs, and at the standardised differences, it is clear that coding was the subtest on which females almost invariably showed better performance. Intriguingly, the female advantage was particularly large in children with mathematic difficulties as compared to children with reading problems or with a mixed profile, with a standardised difference more than twice in the former compared to the latter two groups.

Gender differences in mathematics have been extensively studied. Several studies suggested that boys outperform girls in mathematical problem solving (e.g., Benbow & Stanley, 1983; Geary, 1996; Hyde, Fennema, & Lamon, 1990; for a review, see Hyde, Lindberg, Linn, Ellis, & Williams, 2008), but girls generally outperform boys in arithmetic (e.g., Linn & Hyde, 1989; Willingham & Cole, 1997). Boys' advantage in mathematical problem solving has been attributed to their superior spatial abilities (e.g., Casey, Nuttall, & Pezaris, 1999; Geary, 1996). However, the reasons for girls' advantage in arithmetic are still not clear. One plausible reason for this advantage is that girls have an early advantage in language processing (e.g., Burman, Bitan, & Booth, 2008; Hyde & Linn, 1988). Wei, et al. (2012) tested more than 1500 children aged between 8 and 11 and found that girls outperformed boys in arithmetic tasks. Controlling for scores on the word-rhyming task eliminated gender differences in arithmetic, suggesting that girls' advantage in arithmetic may be related to their advantage in language processing. These results may help us to better understand why, in our sample, differences in coding were higher in children with a specific learning disorder in mathematics than in reading, in which verbal processing is likely to be impaired in both boys and girls.

Despite its theoretical and clinical implications, this study has a number of limitations that should be addressed in future research. Firstly, we only had data on the 10 core subtests of the WISC-IV. It would be particularly interesting to also compare performance on the other additional subtests, despite them not being routinely used by practitioners. Second, we only had data on the WISC-IV; other recently released tests could be used instead, such as the WISC-V (Wechsler, 2014), which is unfortunately still unavailable in Italy. The WISC-V has also introduced some new subtests (e.g., visual puzzles, figure weights, picture span), and it would therefore be interesting to investigate the presence of gender differences in children with SLD on those subtests as well. Finally, we used the WISC, however it would also be interesting to evaluate the progression of gender differences in SLD into adolescents and adults using the WAIS. For these reasons, further study should extend the present research questions to test whether WISC-V processing speed measures can serve as embedded validity indicators (e.g., Erdodi, et al., 2017), also comparing them with other neuropsychological speed tests to further enforce their performance validity. It is also worth mentioning that in our sample the FSIQ was relatively high, as was the vocabulary score. However, vocabulary measures skills that go beyond simple decoding skills, assessing measures of children' verbal fluency, word knowledge, and word usage. Children with SLD are not always impaired on measures of vocabulary as, for example, children tend to also be exposed to words through oral language (Biotteau, et al., 2017; De Clercq-Quaegebeur, et al., 2010; Goswami, et al., 2016; Moura, et al., 2014). Therefore, higher scores on vocabulary and on other measures of the WISC-IV are not entirely surprising. It could also be argued that children with SLD, compared to children with typical development, typically receive higher levels of support and instructions as schooling increases, benefitting their performance on some verbal tasks, which require crystallized intelligence to a larger extent (e.g., vocabulary). However, this hypothesis should be further investigated by future research on this topic.

The identification of children with SLD usually requires assessing intelligence, providing a measure of the FSIQ. In this regard, the new DSM-5 recommendations for the identification of students with SLD advised eliminating the IQ–achievement discrepancy criterion while considering intellectual disability as an exclusion criterion (Tannock, 2013). The IQ-achievement discrepancy is statistically flawed and increases the risk of including students with high IQ but average achievement for their age (see also Snowling, Hulme & Nation, 2020 for the specific case of dyslexia; Mammarella, Toffalini, Caviola, Colling, & Szűcs, 2021, for the specific case of dyscalculia). On the other hand, however, other authors are claiming that if only children with average IQ can be diagnosed with SLD, this will make SLD appear to be a specific disorder (because children with language or attentional problems will tend on average to have lower than average IQs and so fail to be diagnosed as dyslexic) (Snowling, Hulme & Nation, 2020). We believe that the diagnosis of SLD should not rely entirely on the FSIQ but should be based on a clinical evaluation process encompassing a number of indices, including for example the general ability index, which has shown to be less biased in children with various SLDs (see Giofrè et al., 2017 on this point).

We were not particularly interested in differences in terms of statistical significance but in the magnitude of the effect expressed in terms of the effect size. Concerning the analyses on the overall sample, our sample size was sufficient to detect small differences (Bs ≥ .10/.15), but not fully powered to detect very small differences. In any case, smaller differences even if statistically significant would be approaching zero and not of a practical importance in terms of the effect size. While we were interested in understanding differences on factors (e.g., g-factor, PRI, VCI, WMI, and PSI) at the latent level, as for subtests, we were not interested in measuring differences at the latent level, but differences in means at the observed level. Subtests are somewhat less reliable and differences at this level can be attributable to several factors (Styck & Watkins, 2016), and, for this reason, results at the level of the subtests, should be interpreted with care (see Giofrè et al., 2017 on this point). In fact, we found that the intercepts of two subtests (i.e., block design and coding) were not invariant across the two groups. Measurement invariance is determined when a measured construct has the same measurement properties in different groups and, when it is present, we can be certain that any group mean and variance differences in levels of variables marking the construct reflect actual differences in levels of the construct among the groups (Johnson & Bouchard, 2007). When measurement invariance is not achieved, mean and variance differences could reflect differences in the relative importance of the indicators used to measure the construct (Hofer, Horn, & Eber, 1997). Given the body of evidence suggesting that females and males achieve similar levels of overall intellectual processing power using different processes and/or strategies (e.g., Spelke, 2005), it is not unreasonable to hypothesize that measurement invariance across females and males is less than complete (see Johnson & Bouchard, 2007 on this point). Therefore, the lack of invariance on the intercepts of some subtests is not unreasonable and should not be taken as a surprise. As for the analyses on different subgroups, we recognise that those might be important for clinicians, however, dividing the sample into three groups also reduces the statistical power of the analyses. This is an important limitation of the current report, and our results should be replicated using a larger sample of participants further divided in different subgroups of children with SLD.

There are pieces of evidence indicating that gender differences in cognitive abilities vary throughout the life span (Halpern & Wai, 2019). For example, male and female children at younger ages tend to perform equally well on some numerical quantitative tasks while some male advantages emerge later during development particularly when tasks become more complex and tap more heavily into the visuospatial component (Spelke, 2005). It is, however, unclear whether this advantage has a biological origin or is related to the environment. There is evidence indicating that higher exposure to spatially complex environments in males progressively leads to brain changes (Berenbaum et al., 2007). As for children with SLD, to the best of our knowledge, there is no clear data indicating that gender differences tend to vary across the lifespan. This is an interesting hypothesis that should be evaluated by future research, possibly taking a longitudinal approach to evaluate changes in the same participants over time.

To sum up, this paper corroborated the finding that gender differences in general intelligence, if any, are not particularly large. It is most probable that these differences do not reflect real differences in the g-factor but reflect some sort of variation due to the individual subtests included in the battery. Several subtests showed some gender differences, but coding was the subtest with the greatest difference. Investigation of performance on the coding subtest in children with an SLD in mathematics also showed that male children tend to be severely impaired in this subtest, while the performance of females was much higher (almost in the normal range). This, alongside the low *g*-content of this subtest, casts doubt on the inclusion of this subtest in the assessment of general intelligence, which could therefore be biased to some extent.

GENDER DIFFERNCES IN CHILDREN WITH SLD 30

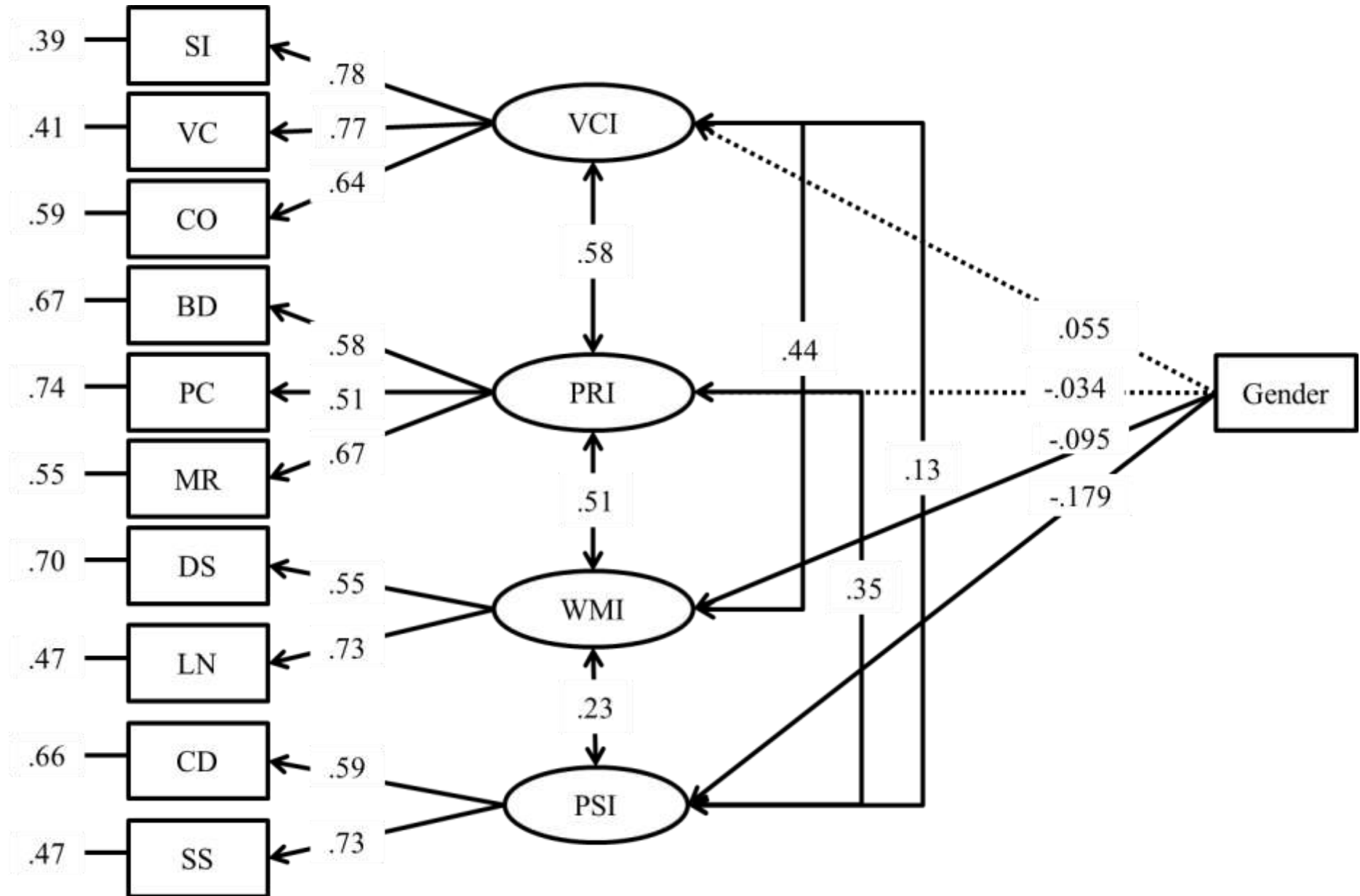

*Figure 1.* Effects of gender on the four principal indices.

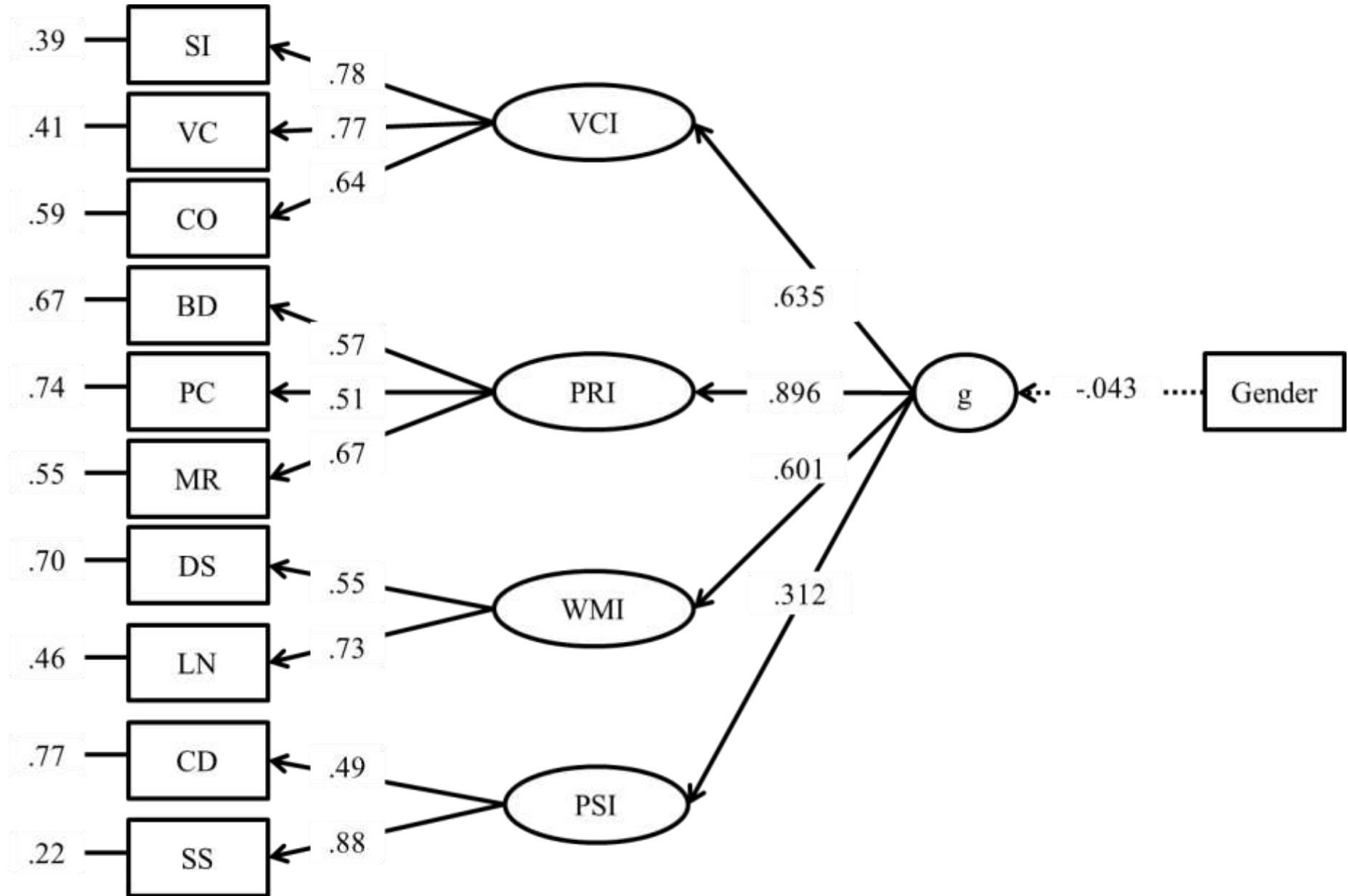

*Figure 2*. Effects of gender on g.

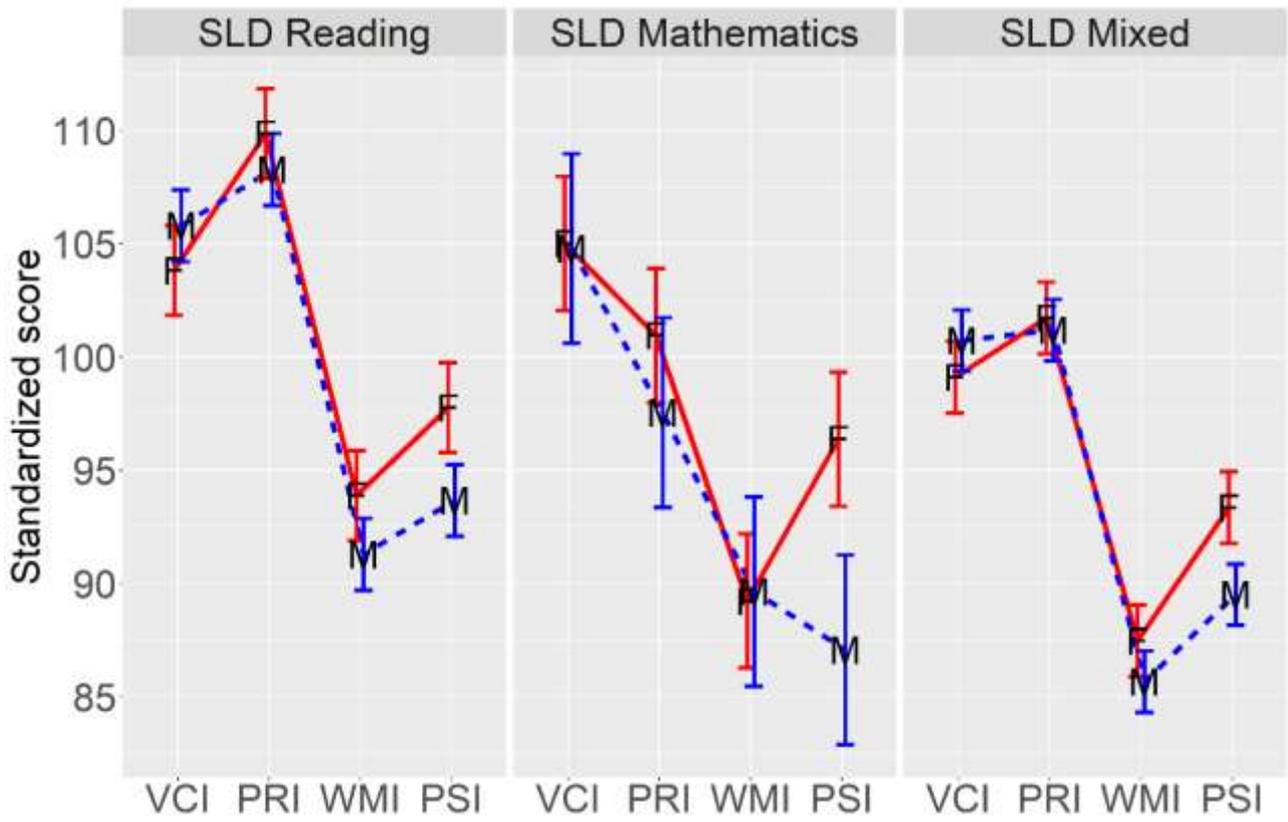

*Figure 3.* Female-male performance on the four principal indices.

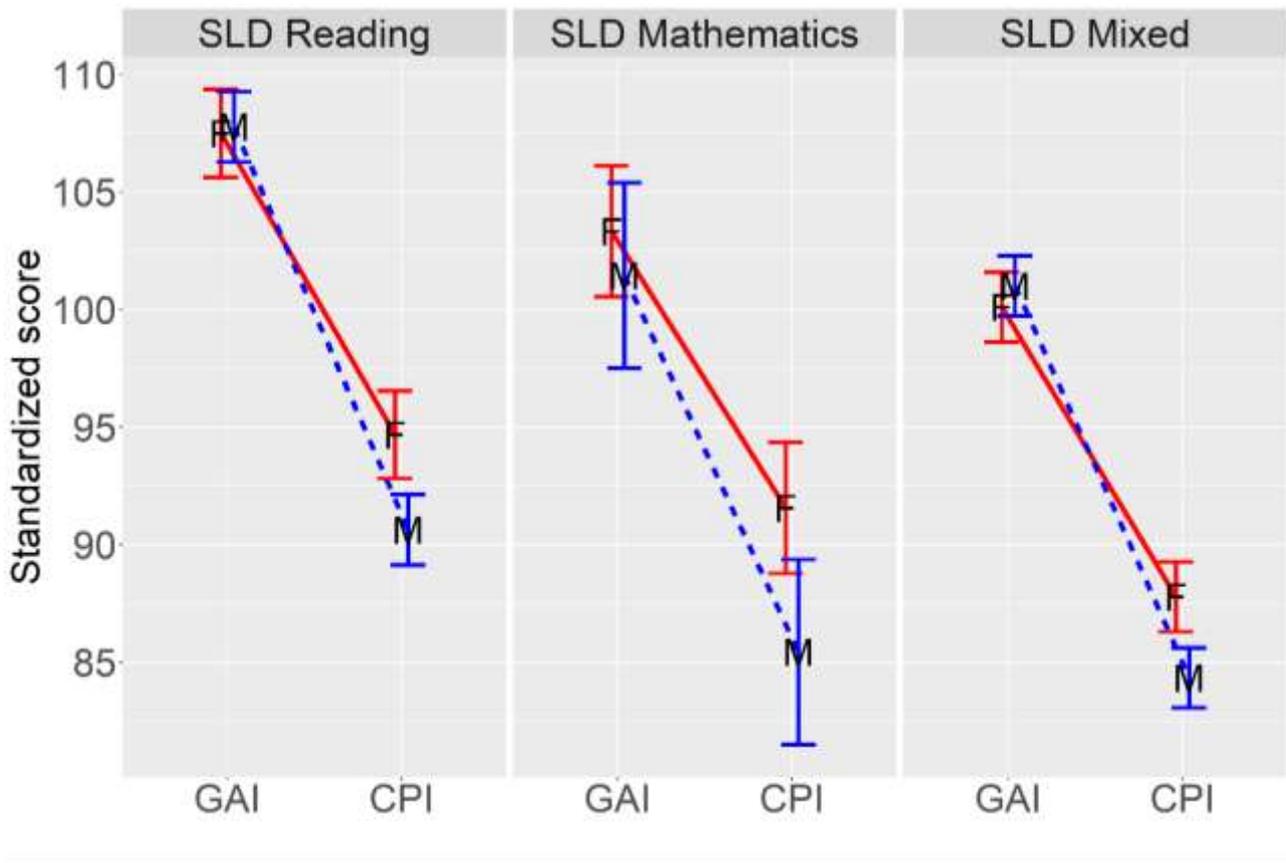

*Figure 4.* Female-male performance on the two ancillary indices.

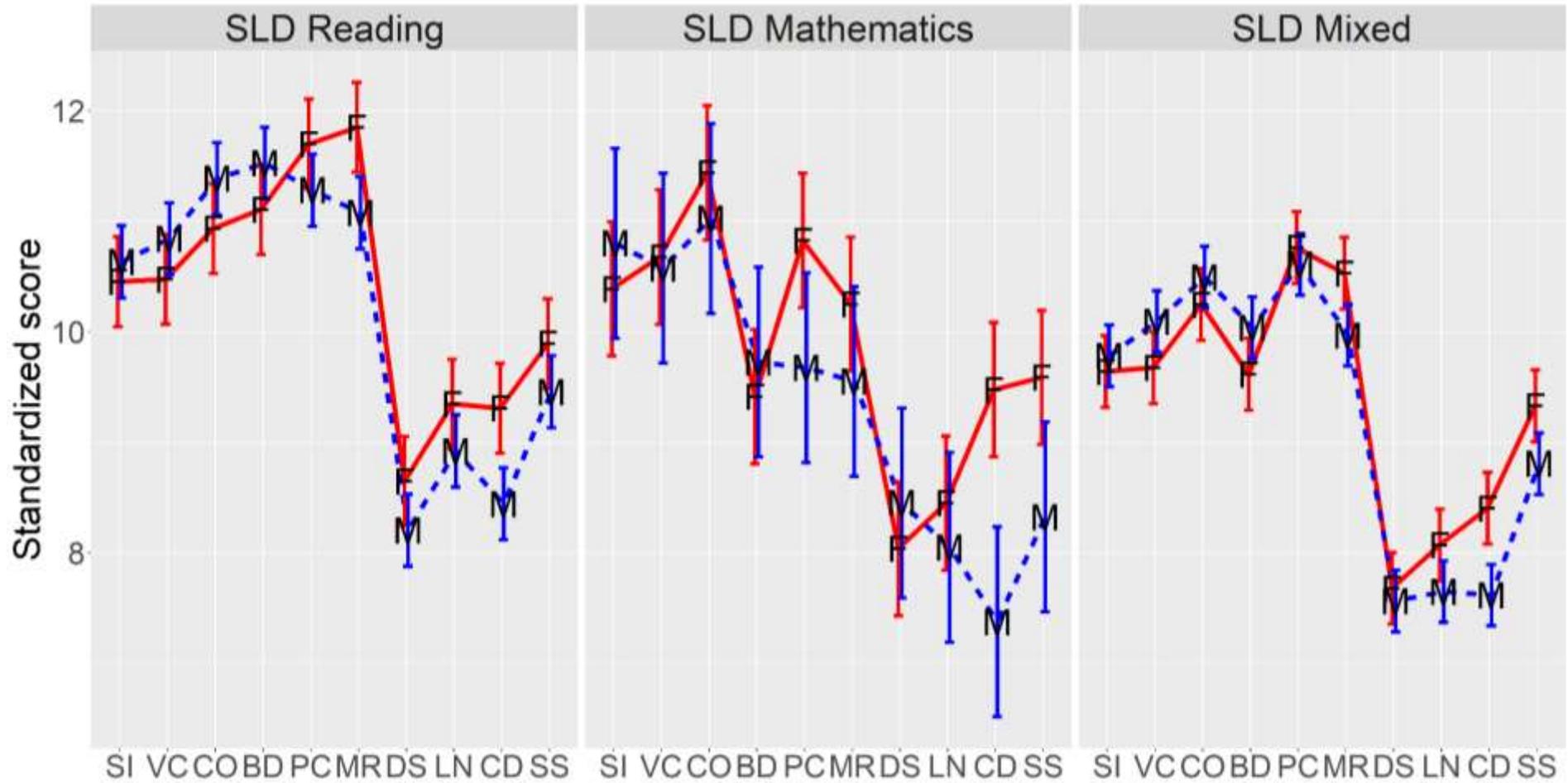

*Figure 5.* Female-male performance on the ten main subtests.

| | | F81.1 Children with dyslexia | | | | | | | F81.2 Children with dyscalculia | | | | | | | F81.3 Children with mixed SLD | | | | | | |
|---|---|---|---|---|---|---|---|---|---|---|---|---|---|---|---|---|---|---|---|---|---|---|
| | | Males | | Females | | | 95% CIs | | Males | | Females | | | 95% CIs | | Males | | Females | | | 95% CIs | |
| | | M | SD | M | SD | d | LL | UL | M | SD | M | SD | d | LL | UL | M | SD | M | SD | d | LL | UL |
| 1 | VCI | 105.81 | 14.08 | 103.84 | 14.96 | 0.136 | -0.056 | 0.330 | 104.8 | 14.9 | 105.03 | 14.11 | -0.016 | -0.405 | 0.375 | 100.73 | 14.91 | 99.11 | 14.64 | 0.110 | -0.044 | 0.263 |
| 2 | PRI | 108.28 | 13.7 | 109.87 | 12.96 | -0.118 | -0.304 | 0.068 | 97.55 | 14.71 | 100.96 | 13.85 | -0.241 | -0.627 | 0.155 | 101.2 | 13.89 | 101.73 | 13.99 | -0.038 | -0.193 | 0.117 |
| 3 | WMI | 91.29 | 11.61 | 93.89 | 13.43 | **-0.211** | -0.401 | -0.015 | 89.65 | 12.04 | 89.24 | 11.33 | 0.036 | -0.352 | 0.433 | 85.67 | 11.79 | 87.48 | 12.17 | -0.151 | -0.306 | 0.005 |
| 4 | PSI | 93.67 | 12.36 | 97.77 | 13.75 | **-0.317** | -0.510 | -0.123 | 87.08 | 15.22 | 96.39 | 15.62 | **-0.601** | -0.990 | -0.195 | 89.52 | 13.11 | 93.35 | 13.8 | **-0.286** | -0.441 | -0.129 |
| 5 | FSIQ | 101.35 | 11.72 | 102.92 | 11.72 | -0.134 | -0.323 | 0.053 | 94.72 | 12.12 | 98.62 | 11.95 | -0.325 | -0.700 | 0.070 | 94.01 | 11.52 | 95 | 11.98 | -0.085 | -0.241 | 0.069 |
| 6 | GAI | 107.77 | 13.39 | 107.5 | 13.12 | 0.021 | -0.170 | 0.207 | 101.45 | 13.06 | 103.33 | 12.87 | -0.146 | -0.528 | 0.243 | 101.01 | 13.38 | 100.11 | 14.53 | 0.065 | -0.093 | 0.219 |
| 7 | ICC | 90.64 | 11.17 | 94.68 | 13.04 | **-0.338** | -0.530 | -0.143 | 85.44 | 13.19 | 91.56 | 12.17 | **-0.489** | -0.885 | -0.082 | 84.34 | 11.41 | 87.78 | 11.65 | **-0.299** | -0.455 | -0.141 |
| 8 | SI | 10.63 | 2.81 | 10.45 | 2.91 | 0.063 | -0.127 | 0.252 | 10.8 | 2.74 | 10.39 | 2.74 | 0.150 | -0.232 | 0.525 | 9.78 | 2.99 | 9.64 | 2.89 | 0.049 | -0.106 | 0.201 |
| 9 | VC | 10.84 | 2.64 | 10.47 | 2.68 | 0.138 | -0.051 | 0.328 | 10.57 | 2.56 | 10.68 | 2.93 | -0.036 | -0.399 | 0.334 | 10.09 | 2.83 | 9.68 | 2.75 | 0.149 | -0.005 | 0.301 |
| 10 | CO | 11.38 | 3.1 | 10.93 | 3.19 | 0.144 | -0.048 | 0.336 | 11.03 | 3.39 | 11.44 | 2.82 | -0.136 | -0.549 | 0.273 | 10.49 | 3.21 | 10.25 | 3.24 | 0.077 | -0.078 | 0.233 |
| 11 | BD | 11.53 | 2.71 | 11.11 | 2.63 | 0.157 | -0.031 | 0.344 | 9.72 | 2.9 | 9.41 | 2.71 | 0.113 | -0.284 | 0.506 | 10.04 | 2.84 | 9.61 | 2.63 | **0.154** | 0.000 | 0.307 |
| 12 | PC | 11.28 | 2.94 | 11.7 | 2.69 | -0.147 | -0.331 | 0.036 | 9.68 | 2.89 | 10.82 | 3.46 | -0.351 | -0.723 | 0.028 | 10.61 | 2.92 | 10.76 | 2.88 | -0.052 | -0.205 | 0.102 |
| 13 | MR | 11.08 | 3.06 | 11.85 | 2.84 | **-0.260** | -0.445 | -0.074 | 9.55 | 3.14 | 10.25 | 2.91 | -0.234 | -0.615 | 0.156 | 9.97 | 2.84 | 10.53 | 3.11 | **-0.189** | -0.347 | -0.031 |
| 14 | DS | 8.21 | 2.42 | 8.65 | 2.72 | -0.174 | -0.369 | 0.020 | 8.45 | 2.51 | 8.04 | 2.07 | 0.186 | -0.234 | 0.596 | 7.56 | 2.36 | 7.68 | 2.42 | -0.048 | -0.202 | 0.108 |
| 15 | LN | 8.92 | 2.23 | 9.35 | 2.46 | -0.183 | -0.374 | 0.008 | 8.05 | 2.35 | 8.45 | 2.34 | -0.171 | -0.547 | 0.216 | 7.65 | 2.31 | 8.07 | 2.41 | **-0.180** | -0.334 | -0.023 |
| 16 | CD | 8.44 | 2.62 | 9.31 | 2.97 | **-0.313** | -0.504 | -0.121 | 7.38 | 2.61 | 9.47 | 3.09 | **-0.714** | -1.087 | -0.303 | 7.62 | 2.61 | 8.4 | 2.78 | **-0.293** | -0.453 | -0.136 |
| 17 | SS | 9.46 | 2.44 | 9.89 | 2.63 | -0.174 | -0.363 | 0.020 | 8.32 | 3.22 | 9.59 | 3.29 | -0.386 | -0.770 | 0.010 | 8.81 | 2.62 | 9.33 | 2.73 | **-0.196** | -0.350 | -0.041 |

Table 1. Male minus female standardized differences on indices and subtests in children with dyslexia, dyscalculia and mixed SLD.

**SUPPLEMENTAL MATERIALS**

**for**

**Decoding gender differences: Intellectual profiles of children with specific learning disabilities**

**Gender Differences in WISC-IV Indices Divided by FSIQ Level**

Participants were divided into 3 subgroups: Low FSIQ level (FSIQ < 90; n = 343, including 138 females and 205 males); Medium FSIQ level (FSIQ in [90, 110]; n = 716, including 319 females and 397 males); High FSIQ level (FSIQ > 110; n = 175, including 79 females and 96 males).

First, we calculated Cohen's d (with their 95% CIs) between males and females in each WISC-IV index. Results are reported in the table below. Negative values indicate gender gaps in favor of males, positive values indicate gender gaps in favor of females.

| Index | HIGH level: Cohen's d | MEDIUM level: Cohen's d | LOW level: Cohen's d |
|---|---|---|---|
| VCI | -0.211 (-0.511, 0.090) | -0.217 (-0.365, -0.069) | -0.159 (-0.376, 0.058) |
| PRI | -0.065 (-0.365, 0.234) | -0.021 (-0.169, 0.126) | 0.006 (-0.210, 0.223) |
| WMI | 0.154 (-0.146, 0.454) | 0.155 (0.007, 0.303) | -0.008 (-0.225, 0.208) |
| PSI | 0.432 (0.129, 0.736) | 0.290 (0.141, 0.438) | 0.270 (0.052, 0.487) |

The results presented in the table suggest that the gender differences are slightly more marked in the Medium and High-level subgroups than in the Low-level subgroup. In addition, with regard to the PSI, the High-level subgroup seem to present a larger gender gap than the other two subgroups. However, none of these differences can be considered statistically reliable due to the large overlap in confidence intervals. CIs are especially large in the Low- and High-level subgroups, due to their smaller size.

In addition, we present a figure with the estimated mean scores of females and males in each WISC-IV index, divided by FSIQ level subgroup. As in the manuscript, the mean scores are estimated via mixed-effects models with index and gender as fixed effects, and participants as random effects (with a random intercept). The figure is presented below. The only notable fact is that PSI appears comparatively more deficient in the High level subgroup (roughly at the same level as WMI), but comparatively less deficient in the Low level subgroup, vis-à-vis the rest of the profile. However, this phenomenon is simply a case of regression towards the mean. Since PSI is the index that is least related to the total FSIQ, it is obviously closer to the population average than the other indexes when subgroups are divided by levels of FSIQ.

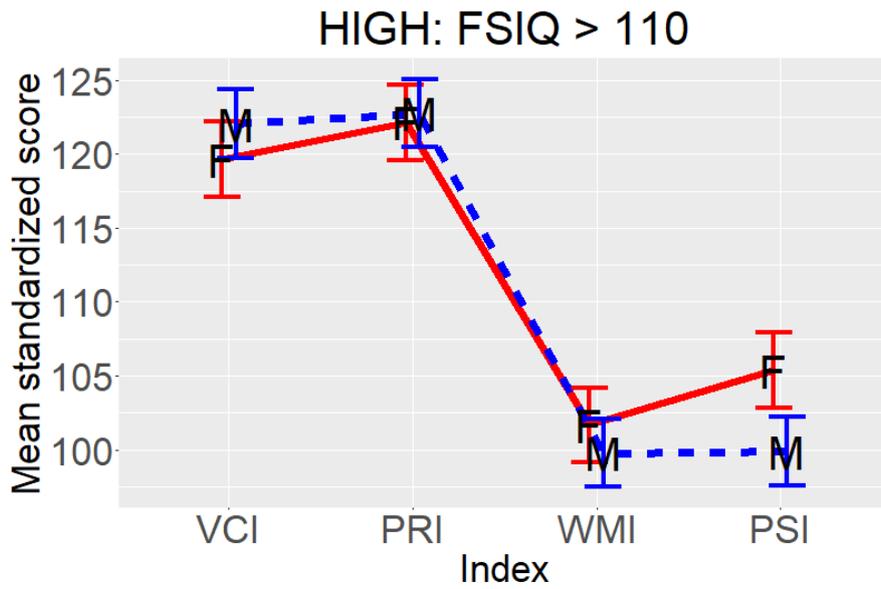

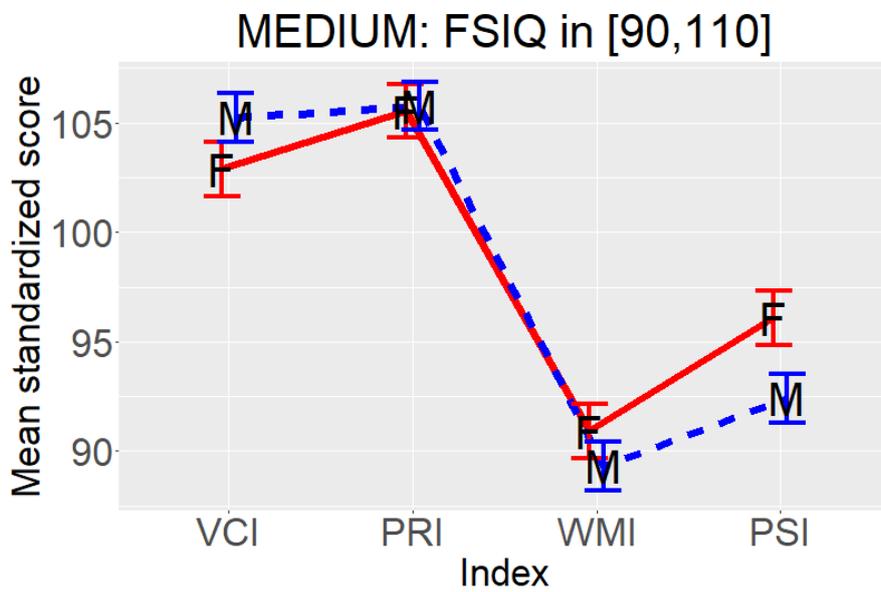

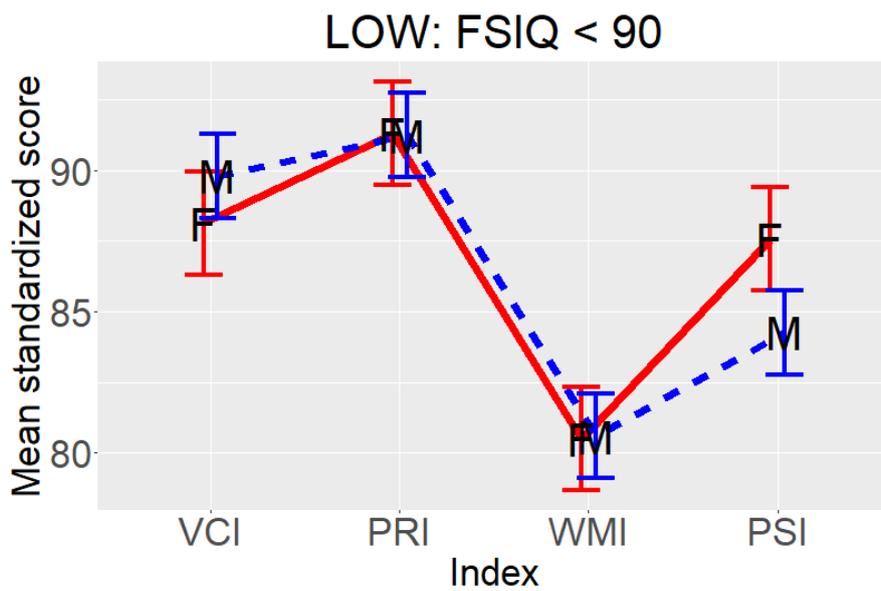